\documentclass[12pt,a4paper]{article}
\usepackage{amsmath}
\usepackage{amssymb}
\usepackage{amsthm}
\usepackage{graphicx}
\usepackage{psfrag}

\newcommand\be{\begin{eqnarray}}
\newcommand\barray{$$\begin{array}{rl}}
\newcommand\ee{\end{eqnarray}}
\newcommand\earray{\end{array}$$}

\newcommand\half{\frac{1}{2}}

\newcommand\inta{\int^{\infty}_{-\infty}}
\newcommand\into{\int^{\infty}_{0}}

\newcommand\dds{\frac{d}{ds}}

\begin{document}
\title{The Basic Elliptic Equations in an Equilateral Triangle}
\author {G. Dassios \\
{\em Division of Applied Mathematics} \\
{\em Department of Chemical Engineering} \\
{\em University of Patras and } \\
{\em ICEHT/FORTH, 26504 Patras, Greece} \\
{\em gdassios@chemeng.upatras.gr} \\
{\em and } \\
{A.S. Fokas} \\
{\em Department of Applied Mathematics} \\
{\em and Theoretical Physics } \\
{\em University of Cambridge } \\
{\em Cambridge, CB30WA, UK} \\
{\em t.fokas@damtp.ac.uk}
\date{May 2004}}

\maketitle

\begin{abstract}
In his deep and prolific investigations of heat diffusion, Lam\'e  was led to the investigation of the
eigenvalues and eigenfunctions of the Laplace operator in an equilateral triangle.  In particular he
derived explicit results for the Dirichlet and Neumann cases using an ingenious change of variables. The
relevant eigenfunctions are complicated infinite series in terms of his variables.

Here we first show that boundary value problems with simple boundary conditions, such as the Dirichlet and
the Neumann problems, can be solved in an elementary manner.  In particular for these problems, the unknown
Neumann and Dirichlet boundary values respectively, can be expressed in terms of a Fourier series. Our
analysis is based on the so called global relation, which is an algebraic equation coupling the Dirichlet
and the Neumann spectral values on the perimeter of the triangle.

As Lam\'e correctly pointed out, infinite series are inadequate for expressing the solution of more complicated
problems, such as mixed boundary value problems. Here we show, utilizing further the global
relation, that such problems can be solved in terms of {\it generalized Fourier integrals}.
\end{abstract}

\section{Introduction}

Solutions of {\it certain} linear elliptic boundary value problems
can be expanded in complete sets of eigenfunctions. Unfortunately,
the actual form of these eigenfunctions is known for only simple
geometries. In fact, only geometries that allow separation of
variables yield well known expressions for the associated
eigenfunctions. But what happens when separation of variables does
not apply? Is it possible to construct the spectral
characteristics of a fundamental domain that does not fit any
separable coordinate system? Some examples where this construction
is possible are presented in the present work. The approach used
here has its roots in the unified transform method for analysing
both lineal and integrable nonlinear PDEs introduced in [4].

A crucial role in this analysis is played by a certain equation
coupling all boundary values, which was called the global relation
in [4]. The concrete form of this equation for the equilateral
triangle was given in the important work of [14], where it was
called a functional equation.

A general overview of the problems solved in this paper is
presented in the sequel where notations and some elementary
formulae are included in order to facilitate the understanding of
the new results. We study boundary value problems for the Laplace,
the Helmholtz and the modified Helmholtz equations in the interior
of an equilateral triangle. These equations are three of the basic
equations of classical mathematical physics. In particular, they
arise as the reduction of several fundamental
parabolic and hyperbolic linear equations. Furthermore, the
specific boundary conditions discussed here cover most cases of
physical significance.

We first introduce some notations.

\subsection{Notations and Useful Identities}

(i) $z$ will denote the usual complex variable and $\alpha$ will denote one of the complex roots of unity,
$$ z = x+iy, \quad \alpha = e^{\frac{2i\pi}{3}} = - \half + \frac{i\sqrt{3}}{2}. \eqno (1.1)$$
Bar will denote complex conjugation, in particular
$$ \bar z = x-iy, \quad \bar \alpha = e^{- \frac{2i\pi}{3}}. $$
$\overline{F(\bar k)}$ will denote the {\it Schwarz conjugate} of the function $F(k)$.

(ii) The complex numbers
$$ z_1 = \frac{l}{\sqrt{3}} e^{\frac{i\pi}{3}}, \quad z_2 = \bar z_1, \quad z_3 = - \frac{l}{\sqrt{3}}, \eqno
(1.2)$$
will denote the vertices of the equilateral triangle, and $D \subset {\mathbb{C}}$ will denote the interior of the
triangle.  The length of each side is $l$.

The sides $(z_2,z_1)$, $(z_3,z_2)$, $(z_1,z_3)$ will be referred to as sides (1), (2), (3) respectively.

%%%%%%%%%%%%%%%%%%%%%%%%%%%%%%%%%%%%%%%%%%%%%%%%%%%%%%%%%%%%%%%%%%%%%%%%%%%%%%%%%%%%%
\begin{center}
\begin{minipage}[b]{6cm}
\psfrag{x}{$x$}
\psfrag{y}{$y$}
\psfrag{a}{$z_{1}$}
\psfrag{b}{$z_{2}$}
\psfrag{c}{$z_{3}$}
\psfrag{A}{$l/2\sqrt{3}$}
\psfrag{B}{$-l/\sqrt{3}$}
\psfrag{T}{$\hat{T}$}
\psfrag{N}{$\hat{N}$}
\centerline{\includegraphics{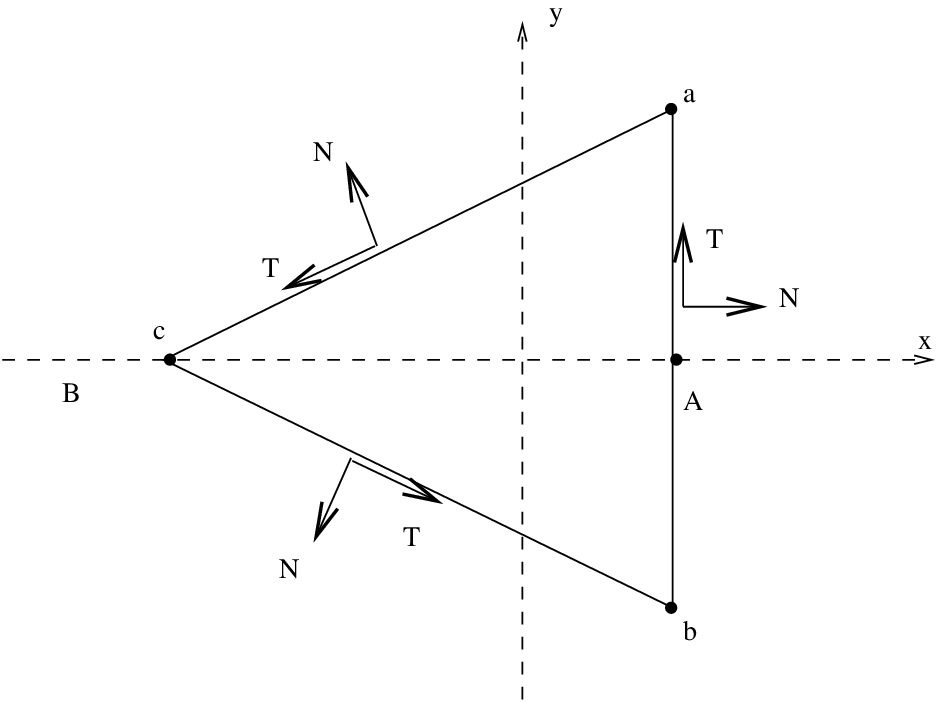}}

\centerline{\textbf{Figure 1.1:} The fundamental domain $D$}
\end{minipage}
\end{center}
%%%%%%%%%%%%%%%%%%%%%%%%%%%%%%%%%%%%%%%%%%%%%%%%%%%%%%%%%%%%%%%%%%%%%%%%%%%%%%%%%%%%%

(iii) On each side we identify the positive direction ${\mathbf{\hat T}}$ and the outward normal $
{\mathbf{\hat N}}$ as in Figure 1.1.
The functions
$$ q^{(j)}(s), \quad q^{(j)}_N(s), \quad s \in \left[ - \frac{l}{2}, \frac{l}{2}\right], \quad j = 1,2,3, \eqno
(1.3)$$
will denote the function $q(x,y)$, as well as the derivative of $q(x,y)$ along the outward normal ${\mathbf{\hat N}}$
 respectively,
for the side $(j)$.

(iv) $E(k)$ and $e(k)$ will denote the following exponential functions
$$ E(k) = \exp \left\{ \left( k + \frac{\lambda}{k} \right) \frac{l}{2\sqrt{3}} \right\},  \quad e(k) =
\exp \left\{ \left( k + \frac{\lambda}{k}\right) \frac{l}{2} \right\}. \eqno (1.4)$$

(v) Using the fact that the numbers $\alpha$ and $\bar \alpha$ satisfy the obvious relations
$$ \alpha^2 = \bar \alpha = \alpha^{-1}, \quad 1+\alpha + \bar \alpha =0, \quad i\bar \alpha - i \alpha =
\sqrt{3}, \quad i\alpha - i = \sqrt{3} \bar \alpha, \eqno (1.5)$$
it is straightforward to obtain
analogous relations for $E(k)$ and $e(k)$.  For example, the last three equations in (1.5)
imply
$$E(k)E(\alpha k)E(\bar \alpha k)=1, \quad E(i\bar \alpha k)E(-i\alpha k)=e(k), \quad E(i\alpha k)E(-ik) = e(\bar
\alpha k). \eqno (1.6)$$

\subsection{Formulation of the Problem}

We will investigate the basic elliptic equations in the interior of the equilateral triangle $D$, namely we will
study the equation
$$q_{xx} + q_{yy} - 4\lambda q =0, \quad (x,y) \in D, \eqno (1.7)$$
where $q(x,y)$ is a real valued function and $\lambda$ is a real constant.  The case of $\lambda =0$, of
$\lambda$ negative, and of $\lambda$ positive, correspond to the Laplace, the Helmholtz, and the modified
Helmholtz equations, respectively.  We will analyze the following problems:

(i) The Dirichlet problem
$$ q^{(j)}(s) = f_j(s), \quad s \in \left[ - \frac{l}{2}, \frac{l}{2} \right], \quad j = 1,2,3 .  \eqno (1.8)$$

(ii) The oblique Robin problem,
$$ \sin \beta q^{(j)}_N(s) + \cos \beta \frac{d}{ds} q^{(j)}(s) + \gamma q^{(j)}(s) = f_j(s), \quad s \in \left[ -
\frac{l}{2}, \frac{l}{2}\right], \quad j = 1,2,3, \eqno (1.9)$$
where $\beta$ and $\gamma$ are real constants and $\sin \beta \neq 0$.  The sum of the first two terms of the
lhs of this equation equals the derivative of $q^{(j)}(s)$ in the direction making an angle $\beta$ with the
positive direction of the side $(j)$.  The Neumann and the Robin problems correspond to the following particular
choices of $\beta$ and $\gamma$,

$$ {\mathrm{Neumann:}} \ \ \beta = \frac{\pi}{2}, \ \ \gamma = 0; \quad {\mathrm{Robin:}} \ \ \beta =
\frac{\pi}{2}, \ \ \gamma \neq 0. \eqno (1.10)$$

(iii) The Poincar\'e type problem
$$\sin\beta_j q^{(j)}_N(s) + \cos \beta_j \frac{d}{ds}q^{(j)}(s) + \gamma_j q^{(j)}(s) = f_j(s), \quad s \in \left[
- \frac{l}{2}, \frac{l}{2}\right], \quad j=1,2,3, \eqno (1.11)$$
where $\beta_1$ is a real constant such that $\sin \beta_1\neq 0$, $\beta_2$ and $\beta_3$ satisfy $\sin\beta_2\neq
0$, $\sin\beta_3 \neq 0$ and are given in terms of $\beta_1$ by the expressions
$$ \beta_2 = \beta_1 + \frac{n\pi}{3}, \quad \beta_3 = \beta_1 + \frac{m\pi}{3}, \quad m,n \in {\mathbb{Z}}, \eqno
(1.12a)$$
and the real constants $\{ \gamma_j\}^3_1$ satisfy the relations
$$ \sin 3\beta_1 \left[ \gamma_2(3\lambda - \gamma^2_2) - e^{in\pi} \gamma_1(3\lambda - \gamma^2_1)\right] =0,
\eqno (1.12b)$$
$$ \sin 3\beta_1 \left[ \gamma_3(3\lambda - \gamma^2_3) - e^{im\pi} \gamma_1(3\lambda - \gamma^2_1)\right] =0.
\eqno (1.12c)$$

A particular case of such a Poincar\'e type problem, which is solved in detail, is the modified Helmholtz equation with Neumann values
on sides (2) and (3) and with Robin values on side (1), where the constant $\gamma$ is given by $\sqrt{3\lambda}$.

We {\it assume} that the functions $f_j$ have sufficient smoothness and that they are compatible at the corners of
the triangle.  The case of boundary conditions which are discontinuous at the corners will be considered elsewhere.

\subsection{The Global Relation}

As it was mentioned earlier the approach used here is based on the analysis of
the global relation, which is the fundamental algebraic relation that couples the Dirichlet and the Neumann values
around the perimeter of the triangle. This equation, first derived for the case of
equilateral triangle in [14] (see also [5]) is
$$E(-ik)\Psi_1(k) + E(-i\bar \alpha k) \Psi_2(\bar \alpha k) + E(-i\alpha k)\Psi_3(\alpha k)   $$
$$=2i \left\{ E(-ik)\Phi_1(k) + E(-i\bar \alpha k)\Phi_2(\bar \alpha k) + E(-i\alpha k)\Phi_3(\alpha k)\right\},
\quad k \in {\mathbb{C}} - \{ 0 \}, \eqno (1.13)$$
where the exponential function $E(k)$ is defined in equation (1.4a), and $\Psi_j$ and $\Phi_j$ are the following
transforms of the Neumann and Dirichlet boundary values:
$$\Psi_j(k) = \int^{\frac{l}{2}}_{- \frac{l}{2}} \exp \left\{ (k+\frac{\lambda}{k})s \right\}
q^{(j)}_N(s)ds, \ \ \Phi_j(k) =
\int^{\frac{l}{2}}_{- \frac{l}{2}} \exp \left\{ (k + \frac{\lambda}{k})s \right\} \left[ \half \dds q^{(j)}(s) + \frac{\lambda}{k}
q^{(j)}(s)\right]ds, \eqno (1.14)$$
for each $j=1,2,3$, and every complex $k\neq 0$.

The general methodology introduced in [4], [5] implies that the global relation must be supplemented by its Schwarz
conjugate, as well as by the four equations obtained from these two equations by replacing $k$ with $\alpha k$ and
with $\bar \alpha k$.  We will refer to these six equations as the {\it basic algebraic relations}.  In this paper
we present two different techniques for solving these equations.

\subsubsection{Solutions via Infinite Series}

For simple problems it is possible to compute the unknown boundary values by evaluating the basic algebraic
relations at particular discrete values of $k$.  This yields the unknown  boundary values in terms of infinite
series.  The Dirichlet and the Neumann problems are examples of problems which
can be solved using this technique.

We use the Dirichlet problem to illustrate this approach.  In this case the functions $\Phi_j$ appearing in the
rhs of the global relation (1.13) can be immediately computed in terms of the given boundary conditions $f_j$, thus
the global relation becomes a single equation for the three unknown functions $\{ \Psi_j\}^3_1$.  Multiplying this
equation by $E(i\alpha k)$, and multiplying the Schwarz conjugate of the global relation by $E(-i\alpha k)$, we
find the following two equations (where we have used the last two of the identities (1.6))
$$e(\bar \alpha k)\Psi_1(k) + e(-k)\Psi_2(\bar \alpha k) + \Psi_3(\alpha k) = 2iA(k), \eqno (1.15)$$
$$e(-\bar \alpha k) \Psi_1(k) + \Psi_2(\alpha k) + e(k) \Psi_3(\bar \alpha k) =-2i B(k). \eqno (1.16)$$
In these equations $A(k)$ and $B(k)$ are known functions and $k\in {\mathbb{C}} - \{ 0 \}$.

For the general Dirichlet problem, we will supplement these two equations with the four equations obtained from
these equations by replacing $k$ with $\alpha k$ and with $\bar \alpha k$.  However, there exists a particular case
for which it is sufficient to analyze only the above two equations.  This is the {\it symmetric} Dirichlet problem,
namely the problem where the functions $f_j$ are all the same, $f_j=f$, $j=1,2,3$.  Then the Neumann values
$q_N^{(j)}(s)$ are also the same, $q_N^{(j)} = q_N$ and hence $\Psi_j(k) = \Psi(k)$, $j=1,2,3$.  Thus
equations(1.15) and (1.16) become two equations for the three unknown functions $\Psi(k)$, $\Psi(\bar \alpha k)$,
$\Psi(\alpha k)$.  Hence, any two of them can be expressed in terms of the remaining one, for example $\Psi(\bar
\alpha k)$ and $\Psi(\alpha k)$ can be expressed in terms of $\Psi(k)$.  In particular, subtracting equations
(1.15), (1.16), we find
$$ \left( e(k)-e(-k)\right) \Psi(\bar \alpha k) = \left( e(\bar \alpha k)-e(-\bar \alpha k)\right) \Psi(k) -
2iG(k), \eqno (1.17)$$
where $G(k) = A(k) + B(k)$ is a known function.  Equation  (1.17) is a single equation for the two unknown
functions $\Psi(\bar \alpha k)$ and $\Psi(k)$.  However, by evaluating this equation at those values of $k$ for
which the coefficient of $\Psi(\bar \alpha k)$ vanishes, i.e. at $e^2(k) =1$, or $k = s_n$,
$$ s_n + \frac{\lambda}{s_n} = \frac{2in\pi}{l}, \quad n \in {\mathbb{Z}}, \eqno (1.18)$$
it follows that $\Psi(s_n)$ can be determined.  Recalling the definition of $\Psi(k)$ and evaluating equation
(1.17) at $k=s_n$, we find
$$\sinh\left[ \left( \bar \alpha s_n + \frac{\lambda}{\bar \alpha s_n}\right) \frac{l}{2} \right]
\int^{\frac{l}{2}}_{- \frac{l}{2}} e^{2i\pi n \frac{s}{l}} q_N(s)ds = iG(s_n), \quad n \in {\mathbb{Z}}. \eqno
(1.19)$$
Thus $q_N(s)$ can be expressed as a Fourier series.

For the general Dirichlet problem (1.8), the six basic algebraic relations couple the nine unknown functions $\{
\Psi_j(k)$, $\Psi_j(\bar \alpha k)$, $\Psi_j(\alpha k)\}^3_1$.  Thus any six of them can be expressed in  terms of
the remaining three.  In particular, it is shown in section 3 that $\Psi_2(\bar \alpha k)$ can be expressed in
terms of $\{ \Psi_j(k)\}^3_1$ by the equation
$$ \left( e^3(-k)-e^3(k)\right) \Psi_2(\bar \alpha k) = \left[ e(-\bar \alpha k) - e^2(-k)e(\bar \alpha k)\right]
\left( \Psi_1(k) + e^2(k) \Psi_3(k)\right) $$
$$ + e^2(-k) \left[ e(-\bar \alpha k) -e^4(k) e(\bar \alpha k)\right] \Psi_2(k) + 2iX(k), \eqno (1.20)$$
where $X(k)$ is known.  In spite of the fact that this equation is a single equation for four unknown functions, it
yields all the three Neumann values $q^{(j)}_N$, $j = 1,2,3$.  Indeed, by evaluating equation (1.20) at those
values of $k$ for which the coefficient of $\Psi_2(\bar \alpha k)$ vanishes, i.e. at $e^6(k) =1$, or $k=k_m$, where

$$ k_m + \frac{\lambda}{k_m} = \frac{2im\pi}{3l}, \quad m \in {\mathbb{Z}}, \eqno (1.21)$$
equation (1.20) yields
$$ \int^{\frac{l}{2}}_{-\frac{l}{2}} e^{\frac{2i\pi m}{3l} s} \left[ q^{(1)}_N(s) + e^{- \frac{2i\pi
m}{3}}q^{(2)}_N(s) + e^{\frac{2i\pi m}{3}} q^{(3)}_N(s)\right] ds = M(k_m), \quad m \in {\mathbb{Z}}, \eqno (1.22)$$
where $M(k_m)$ is known.  This equation in contrast to equation (1.19) involves {\it three} unknown functions.
However, equation (1.21) gives {\it three times} as many values for $m$ as equation (1.18).
Replacing in equation (1.22)
$m$ by $3n$, $3n-1$, $3n-2$, and inverting the left hand sides of the
resulting equations, we find
$$ q_N^{(1)}(s) + q_N^{(2)}(s) + q_N^{(3)}(s) = \frac{1}{l} \sum^\infty_{-\infty} M(3n)e^{- \frac{2i\pi ns}{l}},$$
$$e^{- \frac{2i\pi s}{3l}} \left[ q_N^{(1)}(s) + \alpha q_N^{(2)}(s) + \bar \alpha q_N^{(3)}(s) \right] =
\frac{1}{l} \sum^\infty_{-\infty} M(3n-1) e^{- \frac{2i\pi ns}{l}}, $$
$$ e^{- \frac{4i\pi s}{3l}} \left[ q_N^{(1)}(s) + \bar \alpha q_N^{(2)}(s) + \alpha q_N^{(3)}(s)\right] =
\frac{1}{l}  \sum^\infty_{-\infty} M(3n-2)e^{- \frac{2i\pi ns}{l}}. \eqno (1.23)$$
Thus by solving this system of three algebraic equations it follows that each one of the Neumann boundary values
can be represented in terms of a Fourier series (see Proposition 3.2).

The analysis of the oblique Robin problem (equation (1.9)) is similar.  However, the values of $s_n$ and of $k_m$
in general cannot be found explicitly.  The values $k_m$ satisfy the transcendental equation
$$ e^{\left( k_m + \frac{\lambda}{k_m}\right)l} \frac{ \left( \alpha k_me^{i\beta} + \frac{\lambda}{\alpha k_m
e^{i\beta}} - \gamma \right) \left( \bar \alpha k_me^{-i\beta} + \frac{\lambda}{\bar \alpha k_me^{-i\beta}} -
\gamma \right)}{ \left( \alpha k_me^{-i\beta} + \frac{\lambda}{\alpha k_me^{-i\beta}} - \gamma\right) \left( \bar
\alpha k_me^{i\beta} + \frac{\lambda}{\bar \alpha k_me^{i\beta}}-\gamma\right) } = e^{\frac{2i\pi m}{3}}, \quad m
\in {\mathbb{Z}}. \eqno (1.24)$$
Thus in equation (1.22) instead of $\exp\{ \frac{2i\pi m}{3l}\}$ we now have $\exp\left\{k_m+\frac{\lambda}{k_m}
\right\}$, where $k_m$ satisfies (1.24).
In the particular case of the Neumann problem, $k_m$ satisfies equation (1.21).

\subsubsection{Solutions Via Generalized Fourier Integrals}

For more complicated problems, such as the problem (1.11), the basic algebraic relations can be solved in terms
of a generalized Fourier integral.  This technique in generic, in the sense that it can also be used for the
solution of simple problems.

We use such a simple problem, namely the symmetric Dirichlet problem, to illustrate this approach: It is shown in
Section 4 that the integral defining $\Psi(\bar \alpha k)$ can be solved for $q_N(s)$.  For $\lambda \geq 0$,
$q_N(s)$ is given by
$$q_N(s) = \frac{i\bar \alpha}{4\pi} \int^{\infty e^{\frac{i\pi}{6}}}_{\infty e^{\frac{7i\pi}{6}}}
\exp \left\{-\left( \bar \alpha k + \frac{\lambda}{\bar \alpha k}\right)s \right\}
\left( 1 - \frac{\lambda}{(\bar \alpha k)^2}\right) \Psi(
\bar \alpha k)dk, \quad \lambda \geq 0. \eqno (1.25)$$
Replacing in this equation $\Psi(\bar \alpha k)$ with the expression obtained by solving equation (1.17) for
$\Psi(\bar \alpha k)$, it follows that $q_N(s)$ involves a known integral, as well as an integral containing the
unknown function $\Psi(k)$.  However, using the analyticity properties of the integrant of the latter integral,
it can be shown that this integral can be computed in terms of residues.  Furthermore these residues can be
explicitly calculated in terms of the known function $G(k)$.

The situation for more complicated problems is similar: The unknown boundary values can be expressed in terms of
  known integrals, as well as integrals containing the three unknown functions $\{ \Psi_j(k)\}^3_1$.  Exploiting
the analyticity properties of the integrants of the latter integrals, it can be shown that these integrals can be
computed explicitly.

\subsection{Integral Representations for $q(x,y)$}

When both the Dirichlet and the Neumann boundary values are known, the solution $q(x,y)$ can be determined
either using the classical integral representation in terms of Green's functions [3], or using the novel integral
representations constructed in [5] and [8].  For completeness, both representations are presented in Section 5.

\subsection{Organization of the Paper}

In Section 2 we derive the global relation (1.13).  In Section 3 we solve the symmetric Dirichlet problem
(Proposition 3.1), the general Dirichlet problem (Proposition 3.2),
the general Neumann problem (Proposition 3.3),
and we also discuss the oblique Robin problem. In Section 4 we discuss the basic algebraic relations associated
with the Poincare boundary condition (1.11) and derive the relations (1.12b) and (1.12c).  In Section 5 we
obtain an alternative representation for the symmetric Dirichlet problem and then analyze the
problem defined by equations (1.11) and (1.12).  A particular case of this problem, which is solved in detail
in Proposition 5.1, is a mixed boundary value problem for the modified Helmholtz equation. In Section 6 we discuss the associated
integral representations for $q(x,y)$.  Further discussion of these results is presented in Section 7.

\section{The Global Relation}

Writing the basic elliptic equation (1.7) in the complex variables $(z,\bar z)$ we find
$$q_{z\bar{z}} - \lambda q=0. \eqno (2.1)$$
It is straightforward to verify that this equation can be rewritten in the form [5]
$$ \left( \exp \left\{-ikz - \frac{\lambda}{ik} \bar z \right\} q_z\right)_{\bar{z}} + \frac{\lambda}{ik}
\left( \exp \left\{-ikz - \frac{\lambda}{ik} \bar z \right\} q\right)_{z} = 0, \eqno (2.2)$$
where for the rest of this section $k \in {\mathbb{C}} - \{ 0\}$.  Suppose that equation (2.1) is valid in a simply
connected bounded domain $\Omega \subset {\mathbb{C}}$ with a piecewise smooth boundary $\partial\Omega$.  Then equation
(2.2) and the complex form of Green's theorem imply
$$ \int_{\partial\Omega}
\exp \left\{-ikz - \frac{\lambda}{ik} \bar z \right\}
\left( q_zdz - \frac{\lambda}{ik} qd\bar z\right) =
0. \eqno (2.3)$$
In the particular case that $\Omega$ is the triangular domain $D$, equation (2.3) becomes
$$ \sum^3_{j=1} \tilde\rho_j(k) =0, \eqno (2.4)$$
where the function $\tilde \rho_j(k)$ is given by the following line integral along the side $(j)$ of the
equilateral triangle
$$\tilde \rho_j(k) = \int^{z_j}_{z_{j+1}}
\exp \left\{-ikz - \frac{\lambda}{ik} \bar z\right\}
\left( q_zdz - \frac{\lambda}{ik} qd\bar{z}\right) , \quad j = 1,2,3. \eqno (2.5)$$
In what follows we will show that
$$ \tilde \rho_1(k) = \rho_1(k), \quad \tilde\rho_2(k) = \rho_2(\bar \alpha k), \quad \tilde\rho_3(k) =
\rho_3(\alpha k), \eqno (2.6)$$
where the functions $\rho_j(k)$ are defined in terms of the functions $\Phi_j(k)$ and $\Psi_j(k)$ by the
equation
$$ \rho_j(k) = E(-ik) \left[ \frac{i}{2} \Psi_j(k) + \Phi_j(k)\right], \quad j = 1,2,3. \eqno (2.7)$$

For this purpose we will use the following local parameterizations:

{\it Side 1:} On the side (1) the variable $z$ can be parameterized as
$$ z(s) = \frac{l}{2\sqrt{3}} + is, \quad s \in \left[ -\frac{l}{2}, \frac{l}{2}\right]. \eqno (2.8a)$$
Then $z(-l/2) =z_2$ and $z(l/2) = z_1$.  Since the normal and the tangential derivatives are parallel to the $x$
and to $y$ axes respectively, it follows that
$$ \partial_z = \half \left( \partial_x - i\partial_y\right) = \half \left( \partial_N - i\partial_T\right).
\eqno (2.8b)$$

{\it Side 2:} If $z$ varies along the side (2) and $\zeta$ varies along the side (1), then $z = \zeta\exp\left\{
-i \frac{2\pi}{3}\right\} $.  Thus
$$z(s) = \left( \frac{l}{2\sqrt{3}} + is\right) e^{-i \frac{2\pi}{3}}, \quad s \in \left[ - \frac{l}{2},
\frac{l}{2} \right]. \eqno (2.9a)$$
Note again that $z(-l/2) = z_3$ and $z(l/2) = z_2$.  The equation $\partial_z = \exp\left\{ i \frac{2\pi}{3}\right\}
\partial_\zeta$ implies that
$$ \partial_z = \frac{\alpha}{2} \left( \partial_N - i\partial_T\right). \eqno (2.9b)$$

{\it Side 3:} In analogy with equations (2.9), if $z$ varies along side (3) we find the equations
$$z(s) = \left( \frac{l}{2\sqrt{3}} + is\right) e^{ i \frac{2\pi}{3}}, \quad s \in \left[ - \frac{l}{2},
\frac{l}{2} \right], \eqno (2.10a)$$
and
$$ \partial_z = \frac{\bar\alpha}{2} \left( \partial_N - i\partial_T\right). \eqno (2.10b)$$
Finally, $z(-l/2) = z_1$ and $z(l/2) = z_3$.  Using equations (2.8)-(2.10) in the expressions (2.5) we find
equations (2.6), (2.7).

\section{The Analysis of the Global Relation for Simple \\ Boundary Value Problems}

\subsection{The Symmetric Dirichlet Problem}

We first give the details for the symmetric problem.  In this case
$$ q^{(j)}_N(s) = q_N(s), \quad \Phi_j(k) = F(k), \quad \Psi_j(k) = \Psi(k), \quad j = 1,2,3, \eqno (3.1)$$
where the function $\Psi(k)$ is defined in terms of the unknown function $q_N(s)$ by equation (1.14a) (without
the superscript $(j)$), and the function $F(k)$ is defined in terms of the given boundary condition $f(s)$ by
equation (1.14b), i.e., by the equation
$$ F(k) = \int^{\frac{l}{2}}_{-\frac{l}{2}} \exp \left\{\left( k + \frac{\lambda}{k}\right)s\right\} \left[ \half \dds f(s) +
\frac{\lambda}{k} f(s)\right] ds, \quad k \in {\mathbb{C}} - \{ 0\}. \eqno (3.2)$$
Using equations (3.1), the global relation (1.13) and its Schwarz conjugate yield (1.15) and (1.16), with
$$ A(k) = e(\bar \alpha k)F(k) + e(-k)F(\bar \alpha k) + F(\alpha k),$$
$$ B(k) = e(-\bar \alpha k)F(k) + F(\alpha k) + e(k)F(\bar \alpha k).$$
Hence, since $G=A+B$,
$$ G(s_n) = 2\cosh \left[ \left( \bar \alpha s_n + \frac{\lambda}{\bar \alpha s_n}\right) \frac{l}{2}\right]
F(s_n) + 2e^{in\pi} F(\bar \alpha s_n) + 2F(\alpha s_n). \eqno (3.3)$$
In summary, we have derived the following result:

\paragraph{Proposition 3.1}  Let the real valued function $q(x,y)$ satisfy equation (1.7) in the triangular
domain $D$, with the Dirichlet conditions (1.8), where
$$ f_j(s) = f(s), \quad j=1,2,3, \quad s\in \left[ - \frac{l}{2}, \frac{l}{2} \right], \eqno (3.4)$$
and the function $f(s)$ is sufficiently smooth and satisfies the continuity condition $f(-l/2) = f(l/2)$.  Then
the Neumann boundary values are the same, $q^{(j)}_N(s) = q_N(s)$, $j=1,2,3$, and are given by the Fourier series
$$ q_N(s) = \frac{i}{l} \sum^\infty_{-\infty} e^{- \frac{2in\pi s}{l}} \frac{G(s_n)}{\sinh \left[ \left( \bar
\alpha s_n + \frac{\lambda}{\bar \alpha s_n}\right) \frac{l}{2}\right] }, \eqno (3.5)$$
where $s_n$ is defined by equation (1.18) and $G(s_n)$ is given in terms of $f(s)$ by equations (3.2) and (3.3).

\subsection{The General Dirichlet Problem}

The global relation and its Schwarz conjugate yield equations (1.15) and (1.16), where the known functions $A(k)$
and $B(k)$ are now given by the equations
$$A(k) = e(\bar \alpha k)F_1(k) + e(-k)F_2(\bar \alpha k) + F_3(\alpha k),$$
$$ B(k) = e(-\bar \alpha k) F_1(k) + F_2(\alpha k) + e(k)F_3(\bar \alpha k), \eqno (3.6)$$
where
$$ F_j(k) = \int^{\frac{l}{2}}_{- \frac{l}{2}} \exp \left\{\left( k+\frac{\lambda}{k}\right)s\right\} \left[ \half \dds f_j(s) +
\frac{\lambda}{k} f_j(s)\right] ds, \quad j=1,2,3,\quad k \in {\mathbb{C}} -\{ 0\}. \eqno (3.7)$$
Replacing in equations (1.15) and (1.16) $k$ by $\bar \alpha k$ and then eliminating $\Psi_1(\bar \alpha k)$ from
the resulting two equations we find
$$ e(-\alpha k) \Psi_3(k) + e(k) \Psi_2(\alpha k) - 2ie(-\alpha k)A(\bar \alpha k)   \qquad \ \ \qquad \ \ $$
$$ = e(\alpha k) \Psi_2(k) + e(-k)\Psi_3(\alpha k) + 2ie(\alpha k)B(\bar \alpha k). \eqno (3.8)$$
Taking the Schwarz conjugate of this equation (or equivalently eliminating $\Psi_1(\alpha k)$ from the equations
obtained from equations (1.15) and (1.16) by replacing $k$ with $\alpha k$) we find
$$ e(-\bar \alpha k) \Psi_3(k) + e(k)\Psi_2(\bar \alpha k) + 2ie(-\bar \alpha k) \overline{ A(\bar \alpha \bar
k)}   \qquad \ \ \qquad \ \ $$
$$ = e(\bar \alpha k) \Psi_2(k) + e(-k)\Psi_3(\bar \alpha k) -2ie(\bar \alpha k) \overline{ B(\bar \alpha \bar k)}.
\eqno (3.9)$$
Substituting $\Psi_2(\alpha k)$ from equation (3.8) and $\Psi_3(\bar \alpha k)$ from equation (3.9) into equation
(1.16), we find an equation involving $\Psi_3(\alpha k)$, $\Psi_2(\bar \alpha k)$, and $\{ \Psi_j(k)\}^3_1$.
Eliminating $\Psi_3(\alpha k)$ from this equation and from equation (1.15) we find equation (1.20) with $X(k)$
given by the following equation,
$$ \begin{array}{rl}
X(k) & = \left[ e^2(-k)e(\bar \alpha k) + e(-\bar \alpha k)\right] \left[ F_1(k) + e^2(k) F_3(k)\right] \\ \\
& + e^2(-k)
\left[ e^2(-k)e(\bar \alpha k)e^6(k) + e(-\bar \alpha k)\right] F_2(k)\\ \\
& + 2e^2(k) F_1(\alpha k) + 2F_2(\alpha k) + 2e^2(-k)F_3(\alpha k)\\ \\
& + e^3(-k) \left[ 2e^2(k) F_1(\bar \alpha k) + (e^6(k) + 1)F_2(\bar \alpha k) + 2e^2(-k)e^6(k)F_3(\bar \alpha
k)\right]. \end{array} \eqno (3.10a)$$
Letting $k=k_n$ we find
$$ \begin{array}{rl}
X(k_n) & = \left[ e^2(-k_n)e(\bar \alpha k_n) + e(-\bar \alpha
  k_n)\right] \times \\ \\
& \ \ \ \ \ \  \left[ F_1(k_n) + e^2(-k_n) F_2(k_n) + e^2(k_n)F_3(k_n)\right] \\ \\
& + 2e^2(k_n) \left[ F_1(\alpha k_n) + e^2(-k_n)F_2(\alpha k_n) + e^2(k_n)F_3(\alpha k_n) \right] \\ \\
& + 2e(-k_n) \left[ F_1(\bar \alpha k_n) + e^2(-k_n)F_2(\bar \alpha k_n) + e^2(k_n) F_3(\bar \alpha k_n) \right]
\end{array} \eqno (3.10b)$$
Solving the algebraic equations (1.23) we find the following result:

\paragraph{Proposition 3.2}  Let the real valued function $q(x,y)$ satisfy equation (1.7) in the triangular
domain $D$, with the boundary conditions (1.8), where the given functions $f_j(s)$ have sufficient smoothness and
are continuous at the vertices.  Then the Neumann data $q^{(j)}_N(s)$, $j=1,2,3$ can be expressed in terms of the
given Dirichlet data by the Fourier series
$$q^{(j)}_N(s) = \frac{1}{3l} \sum^\infty_{n=-\infty} \left[ M(k_{3n}) + c^{(j)}_1 e^{\frac{2i\pi s}{3l}
}M(k_{3n-1}) + c^{(j)}_2 e^{\frac{4i\pi s}{3l}} M(k_{3n-2})\right] e^{- \frac{2i\pi ns}{l}}, \eqno (3.11)$$
where $k_m$ is defined by equation (1.21),
$$ c^{(1)}_1 = c^{(1)}_2 =1, \quad c^{(2)}_1 = c^{(3)}_2 = \bar \alpha, \quad c^{(3)}_1 = c^{(2)}_2 = \alpha ,
\eqno (3.12)$$
$$M(k_m) = \frac{2iX(k_m)}{\bar\alpha^ne(\bar \alpha k_m) - e(-\bar \alpha k_m)} \eqno (3.13)$$
and $X(k_m)$ is defined in terms of $f_j(s)$ by equations (3.7) and (3.10).

\subsection{The Oblique Robin Problem}

Suppose that $q(x,y)$ satisfies the Poincar\'e boundary condition (1.11), i.e.
$$q^{(j)}_N(s) = \frac{1}{\sin \beta_j} \left( f^{(j)} - \cos \beta_j \frac{dq^{(j)}}{ds} -
\gamma_jq^{(j)}\right). \eqno (3.14)$$
Substituting this expression in the definition of $\rho_j(k)$, i.e. in the equation (2.7), we obtain
$$ \rho_j(k) = E(-ik) \int^{\frac{l}{2}}_{- \frac{l}{2}} \exp\left\{\left( k + \frac{\lambda}{k}\right)s\right\} \left[
\frac{i}{2} q^{(j)}_N(s) + \half \frac{d}{ds}q^{(j)}(s) + \frac{\lambda}{k} q^{(j)}(s) \right] ds.$$
Integrating by parts we find the following expression for $\rho_j(k)$:
$$\rho_j(k) = iE(-ik) \left[ H_j(k)Y_j(k) + F_j(k) + C_j(k)\right], \quad j = 1,2,3, \eqno (3.15)$$
where the function $H_j(k)$ is defined by
$$ H_j(k) = ke^{i\beta_j} + \frac{\lambda}{ke^{i\beta_j}} - \gamma_j, \eqno (3.16)$$
the function $F_j(k)$ is defined in terms of the given boundary conditions $f_j(s)$ by the equation
$$F_j(k) = \frac{1}{2\sin\beta_j} \int^{\frac{l}{2}}_{-\frac{l}{2}} \exp\left\{\left( k+\frac{\lambda}{k}\right)s\right\}
f_j(s)ds, \eqno (3.17)$$
the function $C_j(k)$ involves the values of $q(x,y) $ at the vertices,
$$C_j(k) = \frac{e^{i\beta_j}}{2\sin\beta_j} \left[ e(-k)q^{(j)} \left( - \frac{l}{2}\right) - e(k)q^{(j)}\left(
\frac{l}{2}\right)\right], \eqno (3.18)$$
and the function $Y_j(k)$ involves the unknown Dirichlet boundary values,
$$ Y_j(k) = \frac{1}{2\sin \beta_j} \int^{\frac{l}{2}}_{-\frac{l}{2}}
\exp \left\{\left( k+\frac{\lambda}{k}\right)s\right\}
q^{(j)}(s)ds. \eqno (3.19)$$
In equations (3.15)-(3.19), $k$ is complex and $k \neq 0$.

In the particular case of the oblique Robin problem (1.9), $\beta_j=\beta$ and $\gamma_j=\gamma$, $j =
1,2,3,$ thus $H_j(k) = H(k)$, where $H(k)$ is defined by equation (3.16) without the subscript $(j)$.
Substituting the expression for $\rho_j(k)$ (with $H_j=H$) in the global relation (2.4) we find
$$E(-ik)\left[ H(k)Y_1(k) + F_1(k)+C_1(k)\right] + E(-i\bar\alpha k) \left[ H(\bar \alpha k)Y_2(\bar \alpha
k) + F_2(\bar \alpha k) + C_2(\bar \alpha k)\right] $$
$$+ E(-i\alpha k) \left[ H(\alpha k)Y_3(\alpha k) + F_3(\alpha k) + C_3(\alpha k)\right] =0. \eqno (3.20)$$
The contribution from the corner terms $C_j$ cancels.  Indeed, this contribution is proportional to the
following expression
$$E(-ik) \left[ e(-k)q^{(1)} \left( - \frac{l}{2} \right) - e(k) q^{(1)} \left( \frac{l}{2}\right)\right] +
E(-i\bar \alpha k) \left[ e(-\bar \alpha k) q^{(2)} \left( - \frac{l}{2} \right) - e(\bar \alpha k) q^{(2)}
\left( \frac{l}{2}\right)\right]$$
$$ +E(-i\alpha k) \left[ e(-\alpha k) q^{(3)} \left( - \frac{l}{2} \right)-e(\alpha k)q^{(3)}
\left( \frac{l}{2}\right) \right]. \eqno (3.21)$$
But, the assumption of continuity at the vertices implies
$$q^{(1)}\left( -\frac{l}{2}\right) = q^{(2)}\left( \frac{l}{2} \right), \quad q^{(2)}\left( - \frac{l}{2}
\right) = q^{(3)} \left( \frac{l}{2} \right), \quad q^{(3)} \left( - \frac{l}{2} \right) = q^{(1)} \left(
\frac{l}{2} \right). \eqno (3.22)$$
Hence the terms $q^{(1)}(-l/2)$ and $q^{(2)}(l/2)$ in the expression (3.21) cancel iff
$$ E(-ik) e(-k) = E(-i\bar\alpha k)e(\bar \alpha k). \eqno (3.23)$$
This equation is indeed valid, and it is the consequence of the identity
$$ \frac{1}{\sqrt{3}} \left( -ik + \frac{\lambda}{-ik}\right) + \left(-k + \frac{\lambda}{-k} \right) =
\frac{1}{\sqrt{3}} \left( -i \bar\alpha k + \frac{\lambda}{-i\bar\alpha k}\right) + \left( \bar\alpha k +
\frac{\lambda}{\bar\alpha k}\right). $$

Using the fact that the corners term cancel, the global relation and its Schwarz conjugate yield (compare
with equations (1.15) and (1.16)) the following equations:
$$ e(\bar\alpha k) H(k)Y_1(k) + e(-k)H(\bar\alpha k) Y_2(\bar \alpha k) + H(\alpha k)Y_3(\alpha k) =-A(k),$$
$$ e(-\bar\alpha k) \overline{H(\bar k)} Y_1(k) + \overline{H(\bar\alpha \bar k)} Y_2(\alpha k) + e(k)
\overline{H(\alpha \bar k)} Y_3(\bar \alpha k) =-B(k), \eqno (3.24)$$
where $A(k)$ and $B(k)$ are defined by equations (3.6) in terms of $F_j$.

Let
$$P(k) = \frac{H(k)}{\overline{H(\bar k)}}. \eqno (3.25)$$
Following precisely the same steps used for the general Dirichlet problem we find the following expression
for $Y_2(\bar\alpha k)$ in terms of $\{ Y_j(k)\}^3_1$:
$$ \left[ e^3(-k) \frac{P^2(\bar \alpha k)}{P^2(\alpha k)} - e^3(k) \frac{P(\alpha k)}{P(\bar \alpha k)}
\right] \frac{\overline{H(\alpha \bar k)}}{\overline{H(\bar k)}} Y_2(\bar \alpha k) = T(k)  $$
$$+ \left[ e(-\bar \alpha k) - e^2(-k) \frac{P(\bar \alpha k)P(k)}{P^2(\alpha k)} e(\bar \alpha k)\right]
\left[ Y_1(k) + e^2(k) \frac{P(\alpha k)}{P(\bar \alpha k)} Y_3(k) \right]  $$
$$+e^2(-k) \frac{P(\bar\alpha k)}{P(\alpha k)} \left[ e(-\bar \alpha k) - e^4(k) \frac{P^3(\alpha
k)}{P^3(\bar \alpha k)} \frac{P(\bar \alpha k)P(k)}{P^2(\alpha k)} e(\bar \alpha k)\right] Y_2(k), \eqno (3.26)$$
where $T(k)$ is given in terms of $F_j(k)$ by the following equation
$$
\overline{H(\bar k)} T(k) =   e(-k) \left[ E^3(-i\alpha k) - E^3(i\alpha k) \frac{P(\bar \alpha
k)}{P^2(\alpha k)} \right] F_1(k) $$
$$  + \!\left[\! E^3(-i\bar \alpha k) \frac{P(\bar \alpha k)}{P(\alpha k)} - E^3(i\bar\alpha k) \frac{1}{P(\bar
\alpha k)} \!\right]\! F_2(k)
+ e(k) \!\left[\! E^3(-i\alpha k) \frac{P(\alpha k)}{P(\bar \alpha k)} - E^3(i\alpha k) \frac{1}{P(\alpha k)}
\!\right]\! F_3(k) $$
$$  + e^2(k) \frac{P(\alpha k)-1}{P(\bar \alpha k)} F_1(\alpha k)
  + \frac{P(\alpha k)-1}{P(\alpha k)} F_2(\alpha k)
+ e^2(-k) \frac{P(\bar \alpha k)(P(\alpha k)-1)}{P^2(\alpha k)} F_3(\alpha k) $$
$$+ e(-k) \frac{P(\bar \alpha k)-1}{P(\alpha k)} F_1(\bar \alpha k)
+ \left( e^3(k) \frac{P(\alpha k)}{P(\bar \alpha k)} - e^3(-k) \frac{P(\bar \alpha k)}{P^2(\alpha k)}
\right) F_2(\bar \alpha k)
    + e(k) \frac{P(\bar \alpha k)-1}{P(\bar \alpha k)} F_3(\bar \alpha k). \eqno (3.27)$$
If $k=k_m$, where $k_m$ is defined by
$$ e^6(k) \frac{P^3(\alpha k)}{P^3(\bar \alpha k)} = 1,$$
then
$$e^2(k) \frac{P(\alpha k)}{P(\bar \alpha k)} = e^{\frac{2i\pi m}{3}}, \quad e^2(-k) \frac{P(\bar\alpha
k)}{P(\alpha k)} = e^{- \frac{2i\pi m}{3}}. $$
Thus evaluating equation (3.26) at $k=k_m$, we find the following expression:
$$\int^{\frac{l}{2}}_{- \frac{l}{2}} \exp\left\{\left( k_m + \frac{\lambda}{k_m}\right)s\right\} \left[ q^{(1)}(s) + e^{-
\frac{2i\pi m}{3}} q^{(2)}(s) + e^{\frac{2i\pi m}{3}} q^{(3)}(s)\right] = G(k_m), \quad m \in {\mathbb{Z}}
\eqno (3.28)$$
where
$$ G(k_m) = \frac{2T(k_m)\sin \beta}{\bar \alpha^m \frac{P(k_m)}{P(\alpha k_m)} e(\bar \alpha
k_m) -e(-\bar \alpha k_m) } . \eqno (3.29)$$

Even if one is able to invert equation (3.28) for the bracket appearing in the integrant of the lhs of equation
(3.28), one will find an expression which will involve an infinite series over the transcendental values of $k_m$. Thus
instead of analyzing the relevant inversion we will analyze the general Robin problem using the generalized
Fourier transform approach (see Section 5).

In the case of the Neumann problem, equations (3.28) and (3.29) simplify and yield the following result.

\paragraph{Proposition 3.3}  Let the real valued function $q(x,y)$ satisfy (1.7) in the triangular domain
$D$, with the Neumann boundary conditions
$$ q^{(j)}_N(s) = f_j(s), \quad s \in \left[ - \frac{l}{2}, \frac{l}{2} \right], \quad j = 1,2,3 \eqno
(3.30)$$
where the functions $f_j(s)$ have sufficient smoothness and are continuous at the vertices of the triangle.
Then the Dirichlet data $q^{(j)}(s)$, $j=1,2,3$ can be expressed in terms of the given Neumann data by the
Fourier series
$$ q^{(j)} (s) = \frac{1}{3l} \sum^\infty_{n=-\infty} \left[ N(k_{3n}) + c^{(j)}_1 e^{\frac{2i\pi s}{3l}}
N(k_{3n-1}) + c_2^{(j)} e^{\frac{4i\pi s}{3l}} N(k_{3n-2})\right] e^{- \frac{2i\pi ns}{l}} \eqno (3.31)$$
where $c_i^{(j)}$ are given by (3.12) and
$$N(k_m) = \frac{2T_N(k_m)}{\bar\alpha^m e(\bar \alpha k_m) -e(-\bar\alpha k_m)}. \eqno (3.32)$$
The known function $T_N(k_m)$ is defined by the equation
$$
-\left( ik+\frac{\lambda}{ik} \right) T_N(k) =   e(-k) \left[ E^3(-i\alpha k) + E^3(i\alpha k) \right]
F_1(k)
  + \left[ E^3(-i\bar \alpha k) + E^3(i\bar \alpha k)\right] F_2(k) $$
$$ + e(k) \left[ E^3(-i\alpha k) + E^3(i\alpha k)\right] F_3(k)
  + 2e^2(k) F_1(\alpha k) + 2F_2(\alpha k) + 2e^2(-k)F_3(\alpha k) $$
  $$+ 2e(-k)F_1(\bar \alpha k) + \left( e^3(k) + e^3(-k)\right) F_2(\bar \alpha k) + 2e(k) F_3(\bar \alpha k),
\eqno (3.33) $$
where $F_j(k)$ is given by (3.17) with $\sin \beta_j = 1$.

\section{Poincar\'e Type Boundary Value Problems}

Suppose that $q(x,y)$ satisfies the Poincar\'e type boundary condition (1.11).  Then substituting the
expression $\rho_j(k)$ from equation (3.15) into the global relation (2.4), we find an equation similar with
equation (3.20), where $H(k)$, $H(\bar\alpha k)$ and $H(\alpha k)$ are replaced by $H_1(k)$, $H_2(\bar
\alpha k)$ and $H_3(\alpha k)$, and $F_j$, $C_j$, $Y_j$ are defined by equations (3.17)-(3.19).  Proceeding
as in Section 3.3, in analogy with equation (3.26), we now find
$$D(k)H_2(\bar \alpha k) Y_2(\bar \alpha k) = \sum^3_{j=1} \Gamma_j(k)H_j(k)Y_j(k) + T(k) + C(k), \eqno
(4.1)$$
where $T(k)$ is defined in terms of the known functions $f_j(s)$, $C(k)$ involves the values of $q$ at the
corners, and $D(k)$, $\{ \Gamma_j(k)\}^3_1$ are defined by the following equations:
$$D(k) = \frac{P_1(\bar \alpha k)}{P_2(\alpha k)P_3(\alpha k)} \left[ e^3(-k) - e^3(k) \frac{P_1(\alpha
k)P_2(\alpha k)P_3(\alpha k)}{P_1(\bar\alpha k)P_2(\bar\alpha k)P_3(\bar\alpha k)}\right], \eqno (4.2)$$
$$\Gamma_1(k) = \frac{1}{P_1(k)} \left[ e(-\bar\alpha k) - e^2(-k)e(\bar\alpha k) \frac{P_1(k)P_1(\bar\alpha
k)}{P_2(\alpha k)P_3(\alpha k)} \right], \eqno (4.3a)$$
$$\Gamma_2(k) = e^2(-k) \frac{P_1(\bar\alpha k)}{P_2(k)P_2(\alpha k)} \left[ e(-\bar\alpha k) -
e^4(k)e(\bar\alpha k) \frac{P_2(k)P_2(\alpha k)}{P_1(\bar\alpha k) P_3(\bar \alpha k)} \right], \eqno
(4.3b)$$
$$ \Gamma_3(k) = e^2(k) \frac{P_1(\alpha k)}{P_3(k)P_3(\bar\alpha k)} \left[ e(-\bar\alpha k) - e^2(-k)
e(\bar\alpha k) \frac{P_3(k)P_3(\bar \alpha k)}{P_1(\alpha k)P_2(\alpha k)} \right], \eqno (4.3c)$$
with
$$P_j(k)  = \frac{H_j(k)}{\overline{H_j(\bar k)}}. \eqno (4.4)$$

In order to be able to solve this problem using a generalized  Fourier integral we require that when $D(k)$
vanishes, then $\Gamma_2(k)$ and $\Gamma_3(k)$ are proportional to $\Gamma_1(k)$.  Actually, $\Gamma_3(k)$
is proportional to $\Gamma_1(k)$ for all complex $k$ provided that
$$P_1(k)P_1(\alpha k) P_1(\bar \alpha k) = P_3(k) P_3(\alpha k)P_3(\bar \alpha k). \eqno (4.5a)$$Equating
the brackets appearing in the definitions of $\Gamma_1(k)$ and $\Gamma_2(k)$, and replacing in the resulting
expression $e^6(k)$ by
$$\frac{P_1(\bar\alpha k)P_2(\bar \alpha k) P_3(\bar\alpha k)}{P_1(\alpha k)P_2(\alpha k)P_3(\alpha k)},$$
it follows that $\Gamma_2(k)$ is proportional to $\Gamma_1(k)$ provided that
$$ P_1(k)P_1(\alpha k)P_1(\bar \alpha k) = P_2(k)P_2(\alpha k)P_2(\bar \alpha k). \eqno (4.5b)$$
Equation (4.5b) is valid if the following two equations are valid:

$$\sin 3(\beta_1-\beta_2) = 0, \quad \gamma_2(3\lambda - \gamma^2_2)\sin 3\beta_1 - \gamma_1(3\lambda -
\gamma^2_1) \sin 3\beta_2 =0. \eqno (4.6)$$
Indeed, in order to simplify equation (4.5b) we first compute the product
\\ $H_1(k)H_1(\alpha k)H_1(\bar\alpha k)$,
$$H_1(k)H_1(\alpha k)H_1(\bar\alpha k) = k^3e^{3i\beta_1} + \frac{\lambda^3}{k^3e^{3i\beta_1}} + 3\lambda
\gamma_1-\gamma^3_1. \eqno (4.7)$$
The function $\overline{H_1(\bar k)}$ can be obtained from $H_1(k)$ by replacing $\beta_1$ with $-\beta_1$,
thus $\overline{H_1(\bar k)}$ $ \overline{H_1(\bar\alpha \bar k)}$ $\overline{H_1(\alpha \bar k)}$ is given by
an expression similar to (4.7) with $\beta_1$ replacing by $-\beta_1$.  Hence equation (4.5b) yields
$$ \frac{k^3e^{3i\beta_1} + \frac{\lambda^3}{k^3e^{3i\beta_1}} + 3\lambda\gamma_1-\gamma^3_1}{
k^3e^{-3i\beta_1} + \frac{\lambda^3}{k^3e^{-3i\beta_1}} + 3\lambda\gamma_1-\gamma^3_1} =
\frac{k^3e^{3i\beta_2} + \frac{\lambda^3}{k^3e^{3i\beta_2}} + 3\lambda\gamma_2-\gamma^3_2}{k^3
e^{-3i\beta_2} + \frac{\lambda^3}{k^3e^{-3i\beta_2}} + 3\lambda\gamma_2-\gamma^3_2}. $$
This equation simplifies to the equation
$$ \left( k^6 - \frac{\lambda^6}{k^6}\right) \sin 3(\beta_1-\beta_2) + \left( k^3 -
\frac{\lambda^3}{k^3}\right) \left[ (3\lambda\gamma_2-\gamma^3_2)\sin 3\beta_1 - (3\lambda
\gamma_1-\gamma^3_1)\sin 3\beta_2\right] = 0, $$
which is valid for all $k$ iff equations (4.6) are valid.

Equation (4.6a) implies $\beta_2 = \beta_1 + n\pi/3$, then $\sin 3\beta_2=\sin 3\beta_1\exp\{in\pi\}$ and we
find equation  (1.12b).  Similarly equation (4.5a) yields equation (1.12c).

\paragraph{The Case that the Corner Terms Cancel} \ \

The definition of the corner terms $C_j(k)$, i.e. equation (3.18), shows that $C_j(k)$ involves
$\exp[i\beta_j]/\sin\beta_j$.  Thus the contribution of the corner terms in the global relation (3.20)
vanishes iff
$$ e^{2i\beta_1} = e^{2i\beta_2} = e^{2i\beta_3}. \eqno (4.8)$$

\paragraph{Example  1. }  \ \
$$ \beta_1=\beta_2 = \beta_3 = \frac{2\pi}{3}, \quad \gamma_j \ \ {\mathrm{arbitrary}}. \eqno (4.9)$$
In this case
$$ H_j(k) = k\alpha + \frac{\lambda}{k\alpha} - \gamma_j, \quad P_j(k) = \frac{k\alpha +
\frac{\lambda}{k\alpha} - \gamma_j}{k\bar\alpha + \frac{\lambda}{k\bar\alpha} - \gamma_j}. \eqno (4.10)$$

\paragraph{Example 2.} \ \
$$\beta_1=\beta_2=\beta_3=\beta, \quad \beta \ \ {\mathrm{arbitrary}}, \quad \gamma_2=\gamma_3=0, \quad
\gamma_1 = (3\lambda)^\half, \quad \lambda >0. \eqno (4.11)$$
In this case
$$  H_1(k)  = ke^{i\beta} + \frac{\lambda}{ke^{i\beta}} - \gamma_1, \quad H_2(k) = H_3(k) = ke^{i\beta}
+ \frac{\lambda}{ke^{i\beta}}. \eqno (4.12)$$
In particular if $\beta=\pi/2$, then
$$H_1(k) = i\left( k-\frac{\lambda}{k}\right) - \gamma_1, \quad H_2(k) = H_3(k) = i\left(k -
\frac{\lambda}{k}\right). \eqno (4.13)$$
Thus
$$P_1(k) = \frac{i\left(k - \frac{\lambda}{k}\right)-\gamma_1}{-i\left( k - \frac{\lambda}{k}\right)-\gamma_1}, \quad P_2(k) =
P_3(k) =-1. \eqno (4.14)$$
Hence
$$D(k) = P_1(\bar \alpha k) \left[ e^3(-k) - e^3(k) \frac{P_1(\alpha k)}{P_1(\bar \alpha k)}\right], \eqno (4.15a)$$
$$\Gamma_1(k) = \frac{1}{P_1(k)} \left[ e(-\bar\alpha k) - e^2(-k)e(\bar \alpha k)P_1(k)P_1(\bar \alpha k)\right], \eqno
(4.15b)$$
$$\Gamma_2(k) = e^2(-k)P_1(\bar \alpha k) \left[ e(-\bar\alpha k) + \frac{e^4(k)e(\bar\alpha k)}{P_1(\bar\alpha k)} \right] \eqno
(4.15c)$$
$$\Gamma_3(k) = e^2(k)P_1(\alpha k) \left[ e(-\bar\alpha k) + \frac{e^2(-k)e(\bar\alpha k)}{P_1(\alpha k)} \right]. \eqno
(4.15d)$$
If $P_1(k)$ is defined by (4.14) with $\gamma^2_1 = 3\lambda$, it can be verified that
$$ P_1(k)P_1(\bar\alpha k) = - \frac{1}{P_1(\alpha k)}.$$
Hence
$$\Gamma_3(k) = - \frac{e^2(k)}{P_1(\bar\alpha k)} \Gamma_1(k), \quad \Gamma_2(k)|_{D(k) =0} = -
\frac{e^2(-k)\Gamma_1(k)}{P_1(\alpha k)}. \eqno (4.16)$$

\paragraph{The Laplace Equation} \ \

In the particular case of the Laplace equation with $\gamma_j =0$, it follows that
$P_j = e^{2i\beta_j}$, i.e. $P_j$ is independent of $k$.

\section{The Analysis of the Global Relation Via Fourier Integrals}

In this section we restrict $\lambda$ to be non-negative.  Slightly more complicated
 formulae can be derived for $\lambda < 0$.

We first derive equation (1.25).  Letting $k = |k|e^{\frac{i\pi}{6}}$, the definition of $\Psi(\bar\alpha k)$ yields
$$\Psi(-i|k|) = \int^{\frac{l}{2}}_{- \frac{l}{2}} \exp \left\{\left( -i|k| + \frac{\lambda}{-i|k|}\right)s\right\}
q_N(s)ds. \eqno (5.1)$$
Suppose that $\lambda >0$.  Letting $t(|k|) = |k| - \lambda/|k|$, it
follows that if $|k| \in (0,+\infty)$, then

\noindent $t \in
(-\infty,\infty)$.  Thus inverting equation (5.1) we find
$$q_N(s) = \frac{1}{2\pi} \inta e^{its} \Psi(-i|k|)dt, \quad s \in \left[ -\frac{l}{2}, \frac{l}{2} \right] , \quad q_N(s) =0
\ \ {\mathrm{elsewhere,}}$$
where $|k|$ (in the argument of $\Psi$) is a function of $t$.
  Rewriting $t$ in terms of $|k|$ we find
$$q_N(s) = \frac{1}{2\pi} \into \exp \left\{i\left( |k| - \frac{\lambda}{|k|}\right)s\right\} \Psi(-i|k|) \left( 1+ \frac{\lambda}{|k|^2}\right)
d|k|.$$
Letting $|k| = ke^{-i\pi/6}$ we obtain
$$ q_N(s) = \frac{1}{2\pi} \int_0^{\infty e^{\frac{i\pi}{6}}}
\exp \left\{-\left( \bar\alpha k + \frac{\lambda}{\bar\alpha k}\right)s \right\}
\Psi(\bar\alpha k) \left( ke^{- \frac{i\pi}{6}} + \frac{\lambda}{ke^{-\frac{i\pi}{6}}}\right) \frac{dk}{k}. \eqno (5.2)$$
The rhs of this equation equals
$$ \frac{1}{2\pi} \int^0_{-\infty e^{\frac{i\pi}{6}}}
\exp \left\{- \left( \bar\alpha \xi + \frac{\lambda}{\bar\alpha \xi}\right) s\right\}
\Psi(\bar\alpha \xi) \left( \xi e^{- \frac{i\pi}{6}} + \frac{\lambda}{\xi e^{- \frac{i\pi}{6}}}\right) \frac{d\xi}{\xi}. \eqno (5.3)$$
Indeed, for the derivation of (5.3) we first observe that the function $\Psi(\bar\alpha k)$ remains invariant under the transformation $k\rightarrow
\bar\alpha \lambda/k$.  Thus making the change of variables $k = \bar\alpha \lambda/\xi$ in the rhs of equation (5.2), and using
$$ke^{- \frac{i\pi}{6}} + \frac{\lambda}{ke^{- \frac{i\pi}{6}}} \rightarrow - \left( \xi e^{- \frac{i\pi}{6}} + \frac{\lambda}{\xi e^{-
\frac{i\pi}{6}}}\right), \quad \frac{dk}{k} =- \frac{d\xi}{\xi}, $$
we find the expression (5.3).  Combining (5.2) and (5.3) we obtain
$$ q_N(s) = \frac{1}{4\pi} \int^{\infty e^{\frac{i\pi}{6}}}_{-\infty e^{\frac{i\pi}{6}}}
\exp \left\{- \left( \bar\alpha k + \frac{\lambda}{\bar\alpha k}\right)s\right\}
\Psi(\bar\alpha k) \left( ke^{- \frac{i\pi}{6}} + \frac{\lambda}{ke^{-
\frac{i\pi}{6}}}\right) \frac{dk}{k}. \eqno (5.4)$$
Using $e^{- \frac{i\pi}{6}} = i\bar\alpha$, this equations becomes
equation (1.25). \\
If $\lambda =0$, we set $k =
t \exp\{ i \frac{7\pi}{6}\}$, $t \in {\mathbb{R}}$, and rewrite (1.14a) as
$$ \Psi(it) = \int^{l/2}_{-l/2} e^{its} q_N(s)ds \eqno (5.5)$$
which is inverted to
$$ q_N(s) = \frac{l}{2\pi} \int^{+\infty}_{-\infty} e^{-its} \Psi(it)dt. \eqno (5.6)$$
Replacing in (5.6) $\{it\}$ with $\{\bar\alpha k\}$, we arrive at
$$q_N(s) = \frac{i\bar\alpha}{2\pi} \int^{\infty e^{  \frac{i\pi}{6}}}_{\infty e^{ \frac{i7\pi}{6}}} e^{- \bar
\alpha ks} \Psi(\bar \alpha k)dk.$$
This equation, in comparison to (1.25) misses a factor of 1/2;
this is due to the linearity of the relevant
transformation in this case.

\subsection{The Symmetric Dirichlet Problem}

Solving equation (1.17) for $\Psi(\bar\alpha k)$ and substituting the resulting expression in equation (1.25) we find
$$q_N(s) = \frac{i\bar\alpha}{4\pi} \int^{\infty e^{\frac{i\pi}{6}}}_{\infty e^{\frac{7i\pi}{6}}}
\frac{\exp \left\{- \left( \bar \alpha
k + \frac{\lambda}{\bar \alpha k}\right)s\right\}}{\Delta (k)}  \left[ (e(\bar \alpha k)-e(-\bar\alpha k)) \Psi(k) - 2iG(k)\right]
\left( 1 - \frac{\lambda}{(\bar\alpha k)^2} \right) dk, \eqno (5.7)$$
where $\Delta(k)$ denotes the coefficient of $\Psi(\bar\alpha k)$ in equation (1.17), i.e.
$$\Delta (k) = e(k) - e(-k).$$
The line $\left( \infty e^{\frac{7i\pi}{6}}, \infty e^{\frac{i\pi}{6}}\right) $ splits the complex $k$-plane into the two half
planes,
$$ {\mathcal{D}}^+ = \left\{ k \in {\mathbb{C}}, \quad \frac{\pi}{6} < \arg k < \frac{7\pi}{6} \right\}, $$
$$ {\mathcal{D}}^- = \left\{ k \in {\mathbb{C}}, \quad \frac{7\pi}{6} < \arg k < \frac{13\pi}{6} \right\}.$$
We observe that
 $$ \exp \left\{- \left( \bar\alpha k + \frac{\lambda}{\bar\alpha k}\right)s\right\} e(\bar\alpha k) \ \ {\mathrm{is \ \ bounded \ \ for}} \ \ k
\in {\mathcal{D}}^-,$$
 $$ \exp \left\{- \left( \bar\alpha k + \frac{\lambda}{\bar\alpha k}\right)s\right\} e(-\bar\alpha k) \ \ {\mathrm{is \ \ bounded \ \ for}} \ \ k
\in {\mathcal{D}}^+. \eqno (5.8)$$
Indeed the exponential of (5.8a) involves $\bar \alpha k\left( \frac{l}{2} - s\right)$, and since $l/2 - s \geq 0$, the
exponential in (5.8a) is bounded for Re$(\bar\alpha k) \leq 0$, i.e. in ${\mathcal{D}}^-$.  Similarly, the exponential of (5.8b)
involves $-\bar\alpha k\left( \frac{l}{2} + s\right)$, which since $l/2 + s \geq 0$, is bounded for Re$(\bar\alpha k) \geq 0$,
i.e. in ${\mathcal{D}}^+$.

We also note that $\Psi(k)/\Delta(k)$ is bounded for all $k \in {\mathbb{C}}$, $k \neq s_n$.  Indeed, for Re $k>0$, $\Delta(k)$ is
dominated by $e(k)$, while for Re $k <0$, $\Delta (k)$ is dominated by $e(-k)$, hence
$$\frac{\Psi(k)}{\Delta(k)} \sim \left\{ \begin{array}{ll}
\Psi(k)e(-k), & {\mathrm{Re}} k>0 \\ \\
-\Psi(k)e(k), & {\mathrm{Re}} k < 0. \end{array} \right. \eqno (5.9)$$
Furthermore $\Psi(k) e(-k)$ involves $k(s-l/2)$ which is bounded for Re $k\geq 0$, while $\Psi(k)e(k)$ involves $k(s+l/2)$ which is
bounded for Re $k\leq 0$ (recall that $-\frac{l}{2} \leq  s \leq  \frac{l}{2}$).

The above considerations imply that the parts of the integral (5.7) containing $e(\bar \alpha k)\Psi(k)$ and $e(-\bar\alpha
k)\Psi(k)$ can be computed by using Cauchy's theorem in ${\mathcal{D}}^-$ and ${\mathcal{D}}^+$ respectively.  The associated residues can
be computed as follows: Let $s^+_n$ and $s^-_n$ denote the subsets of $s_n$ in ${\mathcal{D}}^+$ and ${\mathcal{D}}^-$, respectively.
Evaluating equation (1.17) at $k=s^\pm_n$ we find
$$e(\bar\alpha s^-_n) \Psi(s^-_n) = \frac{2iG(s^-_n)}{-e^2(-\bar\alpha s^-_n) + 1}, \quad -e(-\bar\alpha s^+_n) \Psi(s^+_n) =
\frac{2iG(s^+_n)}{1-e^2(\bar \alpha s^+_n)}.$$
Thus
$$ q_N(s) =- \frac{\bar\alpha}{2\pi}
\int^{\infty e^{\frac{i\pi}{6}}}_{\infty e^{\frac{7i\pi}{6}}}
\exp \left\{- \left( \bar\alpha k + \frac{\lambda}{\bar\alpha k}\right)s\right\} \frac{G(k)}{\Delta(k)} \left[ 1- \frac{\lambda}{(\bar\alpha
k)^2} \right] dk$$
$$   -i\bar\alpha \sum_{s^+_n} \exp \left\{-\left( \bar\alpha s^+_n + \frac{\lambda}{\bar\alpha s^+_n}\right)s\right\}
\frac{G(s^+_n)}{\Delta'(s^+_n)\left[ 1-e^2(\bar\alpha s^+_n)\right]} \left[ 1 - \frac{\lambda}{(\bar\alpha s^+_n)^2} \right] $$
$$ +  i\bar\alpha \sum_{s^-_n} \exp \left\{-\left( \bar\alpha s^-_n + \frac{\lambda}{\bar\alpha s^-_n}\right)s\right\}
\frac{G(s^-_n)}{\Delta'(s^-_n)\left[ 1-e^2(-\bar\alpha s^-_n)\right]} \left[ 1 - \frac{\lambda}{(\bar\alpha s^-_n)^2} \right].
\eqno (5.10) $$

\subsection{The Poincar\'e Problem}

Evaluating equation (4.1) at $k=k_m$, where $k_m$ is a zero of $D(k)$, it follows that the unknown terms $Y_j(k)$ appear in the
form
$$\Gamma_1(k_m) H_1(k_m) \left\{ Y_1(k_m) + e^2(-k_m) \frac{P_1(\bar\alpha k_m)P_1(k_m)}{P_2(k_m)P_2(\bar\alpha k_m)}
\frac{H_2(k_m)}{H_1(k_m)} Y_2(k_m)   \right. $$
$$+ \left. e^2(k_m) \frac{P_1(\alpha k_m)P_1(k_m)}{P_3(k_m)P_3(\bar\alpha k_m)} \frac{H_3(k_m)}{H_1(k_m)} Y_3(k_m) \right\}. $$
The  crucial difference of this general case, as compared with the oblique Robin case (1.9), is the following: Using the
definition of $Y_j(k_m)$ we find that the coefficients of $q^{(2)}(s)$ and of $q^{(3)}(s)$ involve in general $k_m$-dependent
expressions, thus it is not clear how the associated integral can be
inverted.  In contrast, equation (4.1) can be solved using the
approach of section (5.1).  The definition of $Y_2(\bar\alpha k)$, i.e. equation (3.19), and equation(1.25), imply
$$q^{(2)}(s) = \frac{i\bar\alpha\sin\beta_2}{2\pi}
\int^{\infty e^{\frac{i\pi}{6}}}_{\infty e^{\frac{7i\pi}{6}}}
\exp \left\{- \left( \bar\alpha k + \frac{\lambda}{\bar\alpha k}\right)s\right\}
\left(1- \frac{\lambda}{(\bar\alpha k)^2} \right) Y_2(\bar\alpha
k)dk. \eqno (5.11)$$
Solving equation (4.1) for $Y_2(\bar\alpha k)$ and substituting the resulting expression in equation (5.11) we find an integral
involving the three unknown functions $\{ Y_j(k)\}^3_1$. The unknown part of this integral involves the factors (5.8) analyzed
already, as well as factors of the type
$$e^2(-k) \frac{Y_j(k)}{D(k)}, \quad \frac{e^2(k)Y_j(k)}{D(k)}. $$
These terms are bounded for all $k \neq k_m$.  Indeed, ignoring the terms involving $P_j(k)$ we find
$$ e^2(-k) \frac{Y_j(k)}{D(k)} \sim \left\{ \begin{array}{ll}
Y_j(k) e(-k) \cdot e^4(-k), & {\mathrm{Re}} k>0 \\ \\
Y_j(k)e(k), & {\mathrm{Re}} k<0, \end{array} \right.$$
$$ e^2( k) \frac{Y_j(k)}{D(k)} \sim \left\{ \begin{array}{ll}
Y_j(k) e(-k) , & {\mathrm{Re}} k>0 \\ \\
Y_j(k)e(k) \cdot e^4(k), & {\mathrm{Re}} k<0, \end{array} \right.$$
which are identical with the expressions (5.9) except for the occurrence of the factors $e^4(-k)$ and $e^4(k)$ for Re $k>0$ and
Re $k<0$, which are bounded.

The above discussion implies that the integral involving
$$ \frac{e(-\bar\alpha k)}{D(k)} \left[ \frac{H_1(k)}{P_1(k)} Y_1(k) + e^2(-k) \frac{P_1(\bar\alpha k)H_2(k)}{P_2(k)P_2(\alpha
k)} Y_2(k) + e^2(k) \frac{P_1(\alpha k)H_3(k)}{P_3(k)P_3(\bar\alpha k)} Y_3(k)\right], \eqno (5.12)$$
can be computed by using Cauchy's theorem in ${\mathcal{D}}^+$.  Evaluating equation (4.1) at $k^+_m$  it follows that the associated
residue equals
$$   \frac{-[T(k^+_m) + C(k^+_m)]}{1 - \frac{ P_1(k^+_m)P_1(\bar\alpha k^+_m)}{P_2(\alpha k^+_m)P_3(\alpha k^+_m)} e^2(-k^+_m)
e^2(\bar\alpha k^+_m)}. \eqno (5.13)$$
Similarly, the integral involving
$$- \frac{P_1(k)P_1(\bar\alpha k)}{P_2(\alpha k)P_3(\alpha k)} \frac{e(\bar\alpha k)}{D(k)} \left[ \frac{e^2(-
k)}{P_1(k)}H_1(k)  Y_1(k)
+ \frac{e^2(k)P_2(\alpha k)P_3(\alpha k)}{P_1(k)P_1(\bar\alpha k)P_3(\bar\alpha k)} H_2(k) Y_2 (k) \right.$$
$$ \left.  + \frac{P_3(\alpha
k)}{P_1(\bar\alpha k)P_1(\alpha k)} H_3(k)Y_3(k) \right] \eqno (5.14)$$
can be computed by using Cauchy's theorem in ${\mathcal{D}}^-$.  Evaluating equation (4.1) at $k^-_m$, it follows that the associated
residue equals
$$ \frac{[T(k^-_m) + C(k^-_m)]}{ 1 - \frac{P_2(\alpha k^-_m)P_3(\alpha k^-_m)}{P_1(k^-_m)P_1(\bar\alpha k^-_m)} e^2(k^-_m)e^2(-\bar\alpha
k^-_m)}. \eqno (5.15)$$

In what follows we give the details for a mixed Neumann-Robin problem.
\vskip .1in

We will consider Example 2, as it is described by (4.11) with $\beta = \frac{\pi}{2}$.  On side
(1) we assume the
Robin condition
$$ q^{(1)}_N(s) + \sqrt{3\lambda} q^{(1)}(s) = f_1(s), \eqno (5.16)$$
and on sides (2) and (3) we assume the Neumann conditions
$$q^{(2)}_N(s) = f_2(s) \eqno (5.17)$$
and
$$q^{(3)}_N(s) = f_3(s). \eqno (5.18)$$
Then $H_j(k)$, $P_j(k)$, $D(k)$ and $\Gamma_j(k)$, are given by (4.13), (4.14), (4.15a) and (4.15b,c,d) respectively.
Furthermore, $\Gamma_3(k)$ is proportional to $\Gamma_1(k)$ (see equation (4.16a)), while $\Gamma_2(k)$ becomes proportional to
$\Gamma_1(k)$ only on those $k_m$'s for which $D(k)$ vanishes.  These are roots of the transcendental equation
$$ \exp\left\{3\left( k+ \frac{\lambda}{k}\right)\right\} =
\frac{\left( k + \frac{\lambda}{k} - \sqrt{\lambda}\right)\left( k +
\frac{\lambda}{k} - 2\sqrt{\lambda}\right)}{\left( k + \frac{\lambda}{k} + \sqrt{\lambda}\right)\left( k + \frac{\lambda}{k} +
2\sqrt{\lambda}\right)}. \eqno (5.19) $$
For $k=k_m$, equation (4.1), in view of (4.16), implies
$$ \Gamma_1(k_m) \left[ H_1(k_m)Y_1(k_m) - \frac{e^2(-k_m)}{P_1(\alpha k_m)} H_2(k_m)Y_2(k_m) - \frac{e^2(k_m)}{P_1(\bar\alpha
k_m)} H_3(k_m)Y_3(k_m)\right] = -T(k_m). \eqno (5.20)$$
By virtue of (4.15b) and the identity
$$P_1(k_m)P_1(\alpha k_m)P_1(\bar\alpha k_m) = -1 \eqno (5.21)$$
equation (5.20) is written as
$$ \left[ e(-\bar\alpha k_m) + e(\bar\alpha k_m)
\frac{e^2(-k_m)}{P_1(\alpha k_m)} \right] \left[
\frac{H_1(k_m)Y_1(k_m)}{P_1(k_m)} \right.$$
$$+e^2(-k_m)P_1(\bar\alpha
k_m)H_2(k_m)Y_2(k_m)+ e^2(k_m) P_1(\alpha k_m)H_3(k_m)Y_3(k_m)]
=- T(k_m). \eqno (5.22)$$

Since the corner terms $C(k)$ vanish, the representation (5.11) and equation (4.1) yield
$$q^{(2)}(s) = \frac{i\bar\alpha}{2\pi} \int^{\infty e^{i\frac{\pi}{6}}}_{\infty e^{i\frac{7\pi}{6}}}
\exp \left\{- \left(\bar\alpha k + \frac{\lambda}{\bar\alpha k}\right)s\right\}
\left( 1 - \frac{\lambda}{(\bar\alpha k)^2}\right)$$
$$ \frac{1}{H_2(\bar\alpha k)D(k)}
\left[ \sum^3_{j=1} \Gamma_j(k)H_j(k)Y_j(k) + T(k)\right] dk. \eqno (5.23)$$
Utilizing the expression (5.22) we arrive at the following result.

\paragraph{Proposition 5.1.} Let the real valued function $q(x,y)$ satisfy equation (1.7) with $\lambda >0$ in the
triangular domain $D$, with the Robin boundary condition (5.16) on side (1) and the Neumann boundary conditions
(5.17) and (5.18) on sides (2) and (3), where the given functions $f_j(s)$ have sufficient smoothness and are continuous at the
vertices.  Then the Dirichlet value on side (2) is given by
$$q^{(2)}(s) = \frac{i\bar\alpha}{2\pi} \int^{\infty e^{i\frac{\pi}{6}}}_{\infty e^{i\frac{7\pi}{6}}}
\exp\left\{-\left(\bar\alpha k+ \frac{\lambda}{\bar\alpha k}\right)s\right\}
\left( 1 - \frac{\lambda}{(\bar\alpha k)^2}\right)
\frac{T(k)}{H_2(\bar\alpha k)D(k)} dk$$
$$ + \bar\alpha \sum_{k^+_m} \frac{\exp\left\{-\left( \bar\alpha k^+_m + \frac{\lambda}{\bar\alpha k^+_m}\right)s\right\}}{H_2(\bar\alpha
k^+_m)D'(k^+_m)} \left( 1 - \frac{\lambda}{(\bar\alpha k^+_m)^2}\right) \frac{T(k^+_m)}{1+ \frac{E^6(i\alpha
k^+_m
)}{P_1(\alpha k^+_m)}} $$
$$ - \bar\alpha \sum_{k^-_m} \frac{ \exp\left\{- \left( \bar\alpha k^-_m + \frac{\lambda}{\bar\alpha k^-_m}\right)s \right\}}{H_2(\bar\alpha
k^-_m)D'(k^-_n)} \left( 1 - \frac{\lambda}{(\bar\alpha k^-_m)^2} \right) \frac{T(k^-_m)}{1 + P_1(\alpha
k^-_m)E^6(-i\alpha k^-_m)}, \eqno (5.24)$$
where $D'(k^\pm_n)$ denotes the derivative of $ D(k)$ evaluated at $k=k^\pm_n$. The summations are taken over all
$k^+_m \in {\mathcal{D}}^+$ and $k^-_m\in {\mathcal{D}}^-$ respectively,
and $T$, $H_2$, $P_1$ are defined by equations (3.27), (3.16), (4.4) respectively.

There exist similar formulas for the Dirichlet values on sides (1) and (3).

\section{The Integral Representations}
If $\lambda \geq 0$ the classical Green's representation is given by [3]
$$ q({\mathbf{r}}) = \frac{1}{2\pi} \int_{\partial D} \left[ K(2\sqrt{\lambda}|{\mathbf{r}}-{\mathbf{r}}'|)\partial_{n'} q({\mathbf{r}}') - q({\mathbf{r}}')\partial_{n'}
K(2\sqrt{\lambda}|{\mathbf{r}}-{\mathbf{r}}'|)\right] dl({\mathbf{r}}') \eqno (6.1)$$
where the integration is over the boundary $\partial D$ of the triangle in the positive direction, $\partial_{n'}$ denotes the
outward normal derivative on $\partial D$, $dl({\mathbf{r}}')$ is the line element along $\partial D$, and $K(x)$ is the modified Bessel
function of the zeroth order and of the second kind for the modified Helmholtz equation.  For the case of the Helmholtz equation
$K(x)$ is proportional to the Hankel function of the zeroth order and of the first kind, while for the Laplace's equation $K(x)$
is proportional to the logarithm of $x$.

For the Laplace equation, the integral representation constructed in [5] is defined as follows:
$$ \frac{\partial q}{\partial z} = \frac{1}{2\pi} \sum^3_{j=1} \int_{l_j} e^{ikz} \tilde \rho_j(k)dk, \quad z \in D, \eqno (6.2)$$
where the contours $l_j$ are the rays from 0 to $\infty$ specified by the arguments $-\pi/2$, $\pi/6$, $5\pi/6$ respectively, and
the functions $\tilde \rho_j(k)$ are defined by equations (2.6) in terms of $\rho_j(k)$, where the latter functions are defined
by equations (2.7), (1.14) with $\lambda =0$.

%%%%%%%%%%%%%%%%%%%%%%%%%%%%%%%%%%%%%%%%%%%%%%%%%%%%%%%%%%%%%%%%%%%%%%%%%%%%%%%%%%%%%
\begin{center}
\begin{minipage}[b]{6cm}
\psfrag{a}{$\pi/6$}
\psfrag{b}{$k_{R}$}
\psfrag{c}{$k_{I}$}
\psfrag{A}{$l_{1}$}
\psfrag{B}{$l_{2}$}
\psfrag{C}{$l_{3}$}
\centerline{\includegraphics{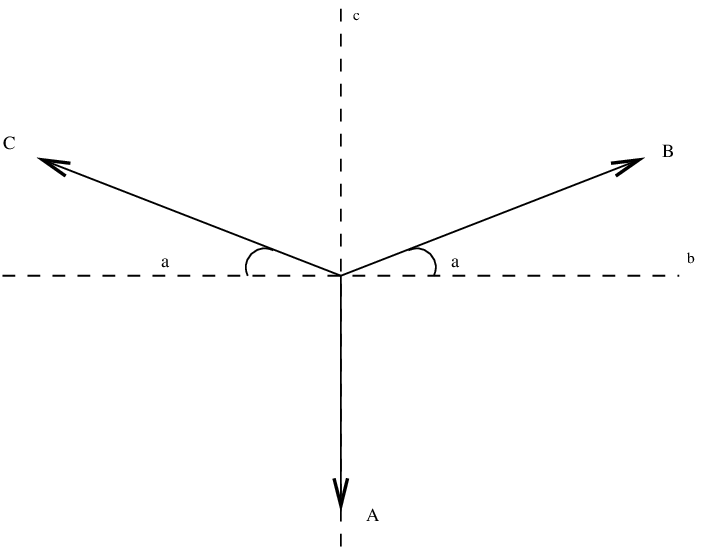}}
\centerline{\textbf{Figure 6.1:} The rays $l_j$ in the complex $k$-plane}
\end{minipage}
\end{center}
%%%%%%%%%%%%%%%%%%%%%%%%%%%%%%%%%%%%%%%%%%%%%%%%%%%%%%%%%%%%%%%%%%%%%%%%%%%%%%%%%%%%%

For the modified Helmholtz equation, the analogue of equation (6.2) is [5]
$$ q(z,\bar z) = \frac{1}{2\pi i} \sum^3_{j=1} \int_{l_j} e^{ikz + \frac{\lambda}{ik}\bar z} \tilde \rho_j(k) \frac{dk}{k}, \quad
z \in D, \eqno (6.3)$$
where the rays $l_j$ are the same as in (6.2) and $\tilde \rho_j(k)$ are defined by equations (2.6) and (2.7).

There exists a similar representation for the Helmholtz equation, which however, in addition to rays, it also involves circular
arcs [5].

\subsection{The Symmetric Dirichlet Problem}
Using the integral representation (6.3) it is possible to compute directly $q(z,\bar z)$, bypassing the computation of the
unknown boundary values.  For brevity of presentation we will only give details for the symmetric Dirichlet problem.
The analysis of the more general boundary value problems (1.8)-(1.11) is similar.

Recalling the definitions of $\tilde \rho_j$, i.e. equations (2.6) and (2.7), it follows that the representation $q(z,\bar z)$
given by equation (6.3) involves the known function $F(k)$ defined by equation (3.2), as well as the unknown
function $\Psi(k)$ which on the rays $l_j$ appears as:
$$ l_1: E(-ik)\Psi(k), \quad l_2: E(-i\bar\alpha k)\Psi(\bar \alpha k), \quad l_3: E(-i\alpha k)\Psi(\alpha k). \eqno (6.4)$$
Solving equation (1.17) for $\Psi(\bar\alpha k)$ in terms of $\Psi(k)$, and then using the Schwarz conjugation of the resulting
equations in order to express $\Psi(\alpha k)$ in terms of $\Psi(k)$, it follows that the expressions in
(6.4) involve the unknown function $\Psi(k)/\Delta(k)$, $\Delta(k)=e(k)-e(-k)$, times the following expressions:
$$l_1: E(-ik) [e(k)-e(-k)], \quad l_2: E(-i\bar\alpha k)[e(\bar\alpha k)-e(-\bar\alpha k)], \quad l_3: E(-i\alpha k)[e(\alpha
k)-e(-\alpha k)]. \eqno (6.5)$$
The third of the relations in (1.5) implies $e(k) = E(i\bar\alpha k)E(-i\alpha k)$, thus
$$ E(-ik)e(k) = E^2(i\bar\alpha k), \quad E(-ik)e(-k) = E^2(i\alpha k).$$
Replacing $k$ by $\bar\alpha k$ and by $\alpha k$ in these identities, we find
$$ E(-i\bar\alpha k)e(\bar\alpha k) = E^2(i\alpha k), \quad E(-i\bar\alpha k)E(-\bar\alpha k) = E^2(ik), $$
$$E(-i\alpha k)e(\alpha k) = E^2(ik), \quad E(-i\alpha k) e(-\alpha k) = E^2(i\bar\alpha k).$$
Thus the expressions in (6.5) involve
$$ l_1: E^2(i\bar\alpha k) - E^2(i\alpha k), \quad l_2: E^2(i\alpha k) - E^2(ik), \quad l_3: E^2(ik) - E^2(i\bar\alpha k).$$
Hence, the unknown part of $q(z,\bar z)$ involves the following integral
$$ J(z,\bar z) = \sum^3_{j=1} J_j (z,\bar z), \eqno (6.6)$$
$$J_3(z,\bar z) = \frac{1}{4\pi} \int_{\{ -l_1\}\cup\{l_2\} } \exp\left\{ik z+
  \frac{\lambda}{i k} \bar z \right\} E^2(i\alpha k)
\frac{\Psi(k)dk}{k\Delta (k)}, \eqno (6.7a)$$
$$J_1(z,\bar z) = \frac{1}{4\pi} \int_{\{ -l_2\}\cup\{l_3\} } \exp \left\{ikz  + \frac{\lambda}{i k} \bar z \right\}
E^2(ik) \frac{\Psi(k)dk}{k\Delta (k)}, \eqno (6.7b)$$
$$J_2(z,\bar z) = \frac{1}{4\pi} \int_{\{ -l_3\}\cup\{l_1\} } \exp \left\{ikz  + \frac{\lambda}{i k} \bar z \right\}
E^2(i\bar \alpha k) \frac{\Psi(k)dk}{k\Delta (k)}. \eqno (6.7c)$$
Each of the above integrals can be computed in terms of residues.  Indeed, it was shown in Section 5 that $\Psi(k)/k\Delta(k)$
is bounded as $k\rightarrow 0$ and as $k\rightarrow \infty$. Furthermore, it will be verified below that the
exponentials,
$$ \exp \left\{ikz + \frac{\lambda}{ik} \bar z\right\} E^2(i\alpha k), \quad
\exp \left\{ikz + \frac{\lambda}{ik}\bar z\right\} E^2(ik), \quad
\exp \left\{ikz + \frac{\lambda}{ik} \bar z\right\} E^2(i\bar \alpha k), \eqno (6.8)$$
are bounded as $k\rightarrow 0$ and $k\rightarrow \infty$, for arg $k$ in
$$ \left[ - \frac{\pi}{2}, \frac{\pi}{6}\right], \quad \left[ \frac{\pi}{6}, \frac{5\pi}{6}\right], \quad \left[ \frac{5\pi}{6},
\frac{3\pi}{2} \right], $$
respectively, provided  that $(z,\bar z) \in D$.  We first consider the first exponential in (6.8); since $z_2 =-\bar\alpha \frac{l}{\sqrt{3}}$,
this exponential can be written as $\exp\{ik(z-z_2) + \lambda(\bar z - \bar z_2)/ik\}$.  If $z$ is in the triangular domain then
$$ \frac{\pi}{2} \leq \arg (z-z_2) \leq \frac{5\pi}{6}.$$
Thus if $- \frac{\pi}{2} \leq \arg k \leq \frac{\pi}{6}$,  we find
$$ 0 \leq \arg [k(z-z_2)] \leq \pi.$$
Hence $\exp\{ik(z-z_2)\}$ is bounded as $|k|\rightarrow\infty$ and $\exp\{\lambda(\bar z - \bar z_2)\bar k/i|k|^2\}$
 is bounded as
$|k|\rightarrow 0$.  Hence there is no contribution from zero and from infinity.

The results for the second and for the third integrals in (6.8) follows
from the above result by using appropriate rotations.
The roots
of $\Delta(k) =0$ lie on the imaginary axis.  Denote by $s^+_n$ those with positive imaginary part and by $s^-_n$ those with
negative imaginary part. Obviously, the residue from each $s^+_n$ has a full contribution to $J_1(z,\bar z)$, while the residue contribution from each
$s^-_n$ is split in two halfs, one half is contributed to $J_2(z,\bar z)$ and one half to $J_3(z,\bar z)$.

Tedious but straightforward calculations lead to the expression
$$J(z,\bar z) = \sum_{s^+_n} \exp\left\{is^+_nz + \frac{\lambda}{is^+_n}\bar z\right\}
\frac{E^2(is^+_n)G(s^+_n)}{s^+_n\Delta'(s^+_n)\Delta(\bar \alpha s^+_n)}$$
$$ + \half \sum_{s^-_n} \exp\left\{is^-_n z + \frac{\lambda}{is^-_n} \bar z\right\} \frac{\left( E^2(i\alpha s^-_n) + E^2(i\bar\alpha
s^-_n)\right)G(s^-_n)}{s^-_n\Delta'(s^-_n)\Delta(\bar\alpha s^-_n)}. \eqno (6.9)$$
In the above calculations, the value of $\Psi_1(s^\pm_n)$ is obtained from (1.18).

We observe that for each $n \in {\mathbb{Z}}$,
$$ \exp \left\{is_nz + \frac{\lambda}{is_n}\bar z\right\} =
\exp \left\{\pm 2\sqrt{(\frac{n\pi}{l})^2+\lambda} x\right\} \cdot
\exp \left\{-2i\frac{n\pi}{l}y\right\} \eqno (6.10)$$
thus this expression shows that the equilateral triangle admits
separable solutions.  It is clear that  each
eigensolution in (6.10),
solves equation (1.7).

\section{Conclusion}
Eigenvalues and eigenfunctions for equation (1.7) with homogeneous Dirichlet, Neumann, and Robin boundary conditions were
constructed in the classical works of Lam\'e [10]-[12].  Some of these results have been rederived by several authors,
in particular the Dirichlet problem is discussed in the recent review [13]. The Robin problem is analysed in [14].
It is remarkable that Lam\'e
  argued, using physical
considerations, that it is impossible  to solve certain problems using
infinite series as opposed to integrals. Indeed Lam\'e writes [12
  p.191]: ``The series should therefore express the fact that the
temperature remains zero on strips of constant width separated by
other strips of double width, in which the temperature may
vary. The analytic interpretation of this sort of discontinuity
demands the introduction of terms where the variables appear {\it inside
integrals}\footnote{These words were not italic in the original
  text}. These terms, of a nature that we will not consider here,
cannot disappear from the total series unless the discontinuity disappears''.

In this paper we have solved several boundary value problems by introducing a novel analysis
of the global relation, i.e. of equation (1.13). Although this equation was first
derived in the important work [14], where it was also used to solve the Robin problem, our treatment of equation
(1.13) is different than that of [14]. As a consequence of our novel analysis of equation (1.13)
we are able to first present a straightforward treatment of {\it simple} boundary value problems.  This treatment,
which is based on the evaluation of the basic algebraic relations (see
the Introduction) at particular values of $k$,
expresses the unknown boundary values in terms of infinite series.
The Dirichlet, Neumann and Robin problems can be solved using this approach.
  We then show that, in agreement with the above remarks of Lam\'e, more {\it complicated}
boundary value problems apparently require the use of generalized Fourier integrals as opposed to infinite series.
Proposition 5.1 presents the solution of such a problem.

In this paper, as opposed to the works of [1], [2], [5]-[7], we have introduced a method
for determining the {\it generalised Dirichlet to Neumann map}, i.e. determining the {\it unknown
boundary values} as opposed to determing $q(x,y)$ itself. In this respect we note that: (a) In some applications one
requires precisely these unknown boundary values.  (b) When both the Dirichlet and the Neumann boundary values are known, it is
straightforward to compute $q(x,y)$.

We emphasize however, that the approach of Section 5 can be used to construct directly
$q(x,y)$.  Indeed, if one uses the novel integral representations for $q(x,y)$ obtained in [5], instead for the representation
(1.25) for $q_N$, and if one follows the approach of Section 5, one can again compute explicitly the contribution of the unknown
functions $Y_j(k)$.  This latter approach is illustrated in Section 6.1 for the symmetric Dirichlet problem.  More complicated
problems using this approach are solved in [2], [6], [7].

In order to compute $q(x,y)$ from the knowledge of both the Dirichlet and the Neumann boundary values one can use either the
classical Green's formulae or the representations of [5].  Regarding the latter representations we note that they provide a
tailor-made transform for the particular problem at hand.  In fact the exponential $\exp\{ikz+(\lambda/ik)\bar z\}$
 reflects the
structure of the PDE, the contours $l_j$ in the complex $k$-plane reflect the geometry of the domain, and the
functions
$\rho_j(k)$ describe the boundary conditions.

Both
the Dirichlet and the Neumann problems involve elementary trigonometric functions. It is interesting that the
analysis of the global condition yields these separable solutions without the direct use of separation of
variables.

For arbitrary values of the constants $\beta_j$ and $\gamma_j$, the Poincar\'e problem (1.11) gives rise to
a matrix Riemann-Hilbert problem.  For the particular case that equations (1.12) are valid, it is possible to
avoid this Riemann-Hilbert problem and to solve the problem in closed form.  Although equations (1.12) impose
severe restrictions on $\beta_j$ and $\gamma_j$, some of the resulting cases appear interesting.  These cases
include the following:

(1) $\beta_1 = \beta_2 = \beta_3 = \frac{2\pi}{3}$, $\gamma_j$ arbitrary.

In this case the angles are specified, but $\gamma_j$ are arbitrary.

(2) The mixed Neumann-Robin problem analyzed in Section 5.

(3) $\beta_2 = \beta   + \frac{4\pi}{3}, \quad \beta_3=\beta + \frac{2\pi}{3}, \quad \gamma_1 = \gamma_2 =
\gamma_3  = \gamma$.

In this case all derivatives are computed along a direction making an angle $\beta$ with the positive vertical
axis.

The results presented  here can be made rigorous, following a formalism similar to the one used in [9].

Several problems remain open which include the following:

1. The investigation of singularities associated with discontinuous boundary conditions.

2. If the  $\beta_j$'s differ, then the global relation (3.20) contains a contribution from $q$ at the three corners.
Several approaches for determining these terms are presented in [1] and [6],  however, the optimal treatment of
these terms remains open.

The approach introduced in [4] and [5] constructs the solutions of a given boundary value problem
{\it without} using eigenfunction expansions. Similar considerations apply to the approach
introduced here for constructing the generalised Dirichlet to Neumann map. However, it turns out that
the above approaches can also be used to investigate the existence of eigenfunction expansions and
to construct these expansions when they exist. This will be presented elsewhere.

\section*{Acknowledgements}

This work was partially supported by the EPSRC. This is part of a joint program undertaken
with A.C. Newell.

\end{document}